\documentclass[11pt]{amsart}
\usepackage{latexsym,amsmath,amsfonts,amssymb, cite,verbatim}
\title[Weyl group multiple Dirichlet series of type $A_{2}$]{Weyl
  group  multiple Dirichlet series of type $A_{2}$}

\author{Gautam Chinta and Paul E. Gunnells}

\address{Department of Mathematics,
The City College of CUNY,
New York, NY 10031, USA
}

\email{chinta@sci.ccny.cuny.edu}

\address{Department of Mathematics and Statistics,
University of Massachusetts,
Amherst, MA 01003
}
\email{gunnells@math.umass.edu}

\date{March 1, 2007}
\thanks{Both authors thank the NSF for support.  The first named
  author gratefully acknowledges the support of the Alexander von
  Humboldt Foundation.}

\DeclareMathOperator{\Supp}{Supp}

\begin{document}

\begin{abstract}
A \emph{Weyl group multiple Dirichlet series} is a Dirichlet series in
several complex variables attached to a root system $\Phi$.  The
number of variables equals the rank $r$ of the root system, and the
series satisfies a group of functional equations isomorphic to the
Weyl group $W$ of $\Phi$.  In this paper we construct a Weyl group multiple
Dirichlet series over the rational function field using $n^{th}$ order
Gauss sums attached to the root system of type $A_{2}$.  The basic
technique is that of \cite{qmds,ratfun}; namely, we construct a
rational function in $r$ variables invariant under a certain action of
$W$, and use this to build a ``local factor'' of the global series.
\end{abstract}

\maketitle

\textit{In memory of Serge Lang}

\newcommand{\Q}{\ensuremath{{\mathbb Q}}}
\newcommand{\F}{\ensuremath{{\mathbb F}}}
\newcommand{\R}{\ensuremath{{\mathbb R}}}
\newcommand{\A}{\ensuremath{{\mathbb A}}}
\newcommand{\C}{\ensuremath{{\mathbb C}}}
\newcommand{\Z}{\ensuremath{{\mathbb Z}}}
\newcommand{\x}{\ensuremath{{\bf { x}}}}
\newcommand{\s}{\ensuremath{{\bf { s}}}}
\newcommand{\p}{\ensuremath{{\mathfrak { p}}}}
\newcommand{\f}{\ensuremath{{\mathfrak { f}}}}
\newcommand{\fS}{\ensuremath{{\mathfrak { S}}}}
\newcommand{\PP}{\ensuremath{{\mathcal P}}}
\newcommand{\QQ}{\ensuremath{{\mathcal Q}}}
\newcommand{\ZZ}{\ensuremath{{\mathcal Z}}}
\newcommand{\TT}{\ensuremath{{\mathcal T}}}

\newcommand{\im}{\mbox{Im}}
\newcommand{\sgn}{\mbox{sgn}}
\newcommand{\res}{\mathop{\mbox{Res}}}
\newcommand{\re}{\mbox{Re}}
\newcommand{\half}{{\smfrac{1}{2}}}

\newcommand{\I}{{\mathcal I}}
\newcommand{\J}{{\mathcal J}}
\newcommand{\N}{{\mathcal N}}
\newcommand{\pr}{{\prime}}
\newcommand{\norm}[1]{\left|\left|#1\right|\right|}
\newcommand{\abs}[1]{\left|#1\right|}

\newcommand{\prs}[2]{\Bigl(\frac{#1}{#2}\Bigr)}
\newcommand{\smfrac}[2]{{\textstyle \frac{#1}{#2}}} 
\newcommand{\nn}{\nonumber}

\newcommand{\OO}{{\mathcal O}}
\newcommand{\OOmon}{{\OO_{\text{mon}}}}
\newcommand{\hpl}{h^{(\p,l)}}
\newcommand{\fpl}{f^{(\p,l)}}

\newtheorem{theorem}{Theorem}
\newtheorem{lemma}[theorem]{Lemma}
\newtheorem{prp}[theorem]{Proposition}
\newtheorem{corollary}[theorem]{Corollary}
\newtheorem{conjecture}[theorem]{Conjecture}
\newtheorem{fact}[theorem]{Fact}
\newtheorem{claim}[theorem]{Claim}
\theoremstyle{definition}
\newtheorem{example}[theorem]{Example}
\newtheorem{remark}[theorem]{Remark}

\newtheorem{defin}[theorem]{Definition}

\numberwithin{theorem}{section}
\numberwithin{equation}{section}

\section{Introduction}
\label{sec:intro}

Weyl group multiple Dirichlet series are Dirichlet series in $r$
complex variables $s_1, s_2, \ldots, s_r$ that have analytic
continuation to $s_i\in \C^r$, satisfy a group of functional equations
isomorphic to the Weyl group of a finite root system of rank $r$, and
whose coefficients are products of $n^{th}$ order Gauss sums.  The
study of these series was introduced in \cite{wmd1}, which also
suggested a method for proving the analytic continuation and 
functional equations of these series.  However, complete
proofs have only been given in certain cases: 
\begin{enumerate}
\item  the \emph{stable
case}, when the rank $r$ of $\Phi$ is small relative to the order $n$
\cite{wmd2}, and 
\item  the \emph{quadratic case}, when $n=2$ and with
$\Phi$ of any rank \cite{qmds, ratfun}.
\end{enumerate} 

In this paper we extend the quadratic case methods of
\cite{qmds,ratfun} to arbitrary $n$ for the root system $A_2$. The
resulting series is not new: for $A_{2}$ and $n\geq 2$ one falls in the stable
range, and therefore our result follows from the work of \cite{wmd2}.
(In fact, this case was treated earlier in \cite{wmd1}.)
Nevertheless there are several reasons why a new treatment of $A_{2}$ is
desirable.  First, the methods used here are completely different from
those of \cite{wmd2} and give a different technique to construct Weyl group
multiple Dirichlet series.  Second, the technique presented here will
work for a root system $\Phi$ of arbitrary rank and for arbitrary $n$,
with no stability restriction.  This will be the subject of a
forthcoming publication; one of the main goals of the present paper is
an exposition of our method in the simplest nontrivial case, namely
$\Phi = A_{2}$. 

With this latter goal in mind we also adopt certain assumptions to make the
exposition simpler.  For instance, we work over a rational function
field to avoid the annoyance of having to deal with Hilbert symbols.
We also focus on the \emph{untwisted case} (see \S \ref{sec:Waction}
for an explanation of this terminology) to avoid some notational
complexities.  A comparison with the methods of \cite{wmd1, qmds,
ratfun} indicates how to extend our methods to an arbitrary
global field containing the $2n^{th}$ roots of unity and to arbitrary
twists.

We now describe our main result in greater detail.  Let $\F$ be a
finite field whose cardinality $q$ is congruent to $1\bmod 4n$. Let $K$
be the rational function field $\F(t)$, and let 
$\OO=\F[t]$.   Let $\OOmon \subset \OO$ be the subset of monic
polynomials.  We let $K_\infty=\F((t^{-1}))$ denote the field of
Laurent series in $t^{-1}$.  

For $x,y\in \OO$ relatively prime, we denote by $\prs {x}{y}$ the
$n^{th}$ order power residue symbol.  We have the reciprocity law
\begin{equation}
  \prs{x}{y}=\prs{y}{x}
\end{equation}
for $x,y$ monic.  The reciprocity law takes this particularly simple
form because of our assumption that the cardinality of $\F$ is
congruent to $1 \bmod 4$.

Let $y\mapsto e(y)$ be an additive character on $K_\infty$ with the
following property: if $I\subset K$ is the set of all $y\in K$
such that the restriction of $e$ to $y\OO$ is trivial, then $I = \OO$.
%
%
Fix an embedding $\epsilon$ from the the $n^{th}$ roots of unity in
$\F$ to $\C^\times.$ For $r,c\in \OO$ we define the Gauss sum $g (r,
\epsilon , c)$ by 
$$g(r,\epsilon, c)=\sum_{y \bmod c} \epsilon\left(\prs{y}{c}\right)
e\left(\frac {ry}{c}\right).$$  
We will also use the notation
$g_i(r, c)=g(r,\epsilon^i, c)$ and $g(r, c)=g(r,\epsilon, c).$
Note that $\epsilon^i$ is not necessarily an embedding.

%

We are now ready to define our double Dirichlet series.  Put 
\begin{multline}\label{def:Z}
Z(s_1,s_2)= \\
(1-q^{n-ns_1})^{-1} (1-q^{n-ns_2})^{-1}
(1-q^{2n-ns_1-ns_2})^{-1}\sum_{\substack {c_1\in \OOmon}} 
\sum_{\substack{c_2\in\OOmon}}  \frac{H(c_1,c_2)}
{|c_1|^{s_1}|c_2|^{s_2}},
\end{multline}
where the coefficient $H(c_1,c_2)$ is defined as follows:
\begin{enumerate}
\item (Twisted multiplicativity) If $gcd(c_1c_2, d_1d_2)=1$ then
\begin{equation}\label{eqn:Htm}
\frac {H(c_1d_1, c_2d_2)}{H(c_1,c_2)H(d_1, d_2)} =
\prs{c_1}{d_1} \prs{d_1}{c_1} \prs{c_2}{d_2} \prs{d_2}{c_2}
\prs{c_1}{d_2} ^{-1}\prs{d_1}{c_2} ^{-1}.
\end{equation}
\item ($\p$-part) If $\p$ is prime, then
\begin{multline}\label{eqn:Hp}
\sum_{k,l\geq 0}H(\p^k, \p^l)x^ky^l=
1+g(1, \p) x+ g(1, \p) y+
g(1, \p) g(\p, \p^2)xy^2\\+
g(1,\p)g(\p, \p^2) x^2 y+ g(1,  \p)^2g(\p, \p^2) x^2 y^2.
\end{multline}
\end{enumerate}

Our main result is

\begin{theorem}\label{thm:main}
  The double Dirichlet series $Z(s_1, s_2)$ converges absolutely for
$\re(s_i)$ sufficiently large and has an analytic continuation to
all $(s_1,s_2)\in \C^{2}$. Moreover, $Z(s_1, s_2)$ satisfies two functional
equations of the form 
\begin{equation}\label{eqn:sigmas} \sigma_{1}\colon (s_1, s_2)\mapsto (2-s_1,
s_1+s_2-1)\text{\ and\ } \sigma_{2}\colon (s_1, s_2)\mapsto ( s_1+s_2-1, 2-s_2).
\end{equation}
These two functional equations generate a subgroup of the affine
transformations of $\C^{2}$ isomorphic to the symmetric group $S_3$.
\end{theorem}

The precise statement of the functional equations involves a set of
double Dirichlet series $Z (s_{1},s_{2};i,j)$, where $0\leq i,j\leq
n-1$, and where $Z (s_{1},s_{2}) = \sum_{i,j} Z (s_{1},s_{2};i,j)$; we
refer to Theorem \ref{thm:final} for details.  Moreover, one can
explicitly write down $Z(s_1,s_2)$ as a rational function in
$q^{-s_1}, q^{-s_2}.$ For $n=2,$ this was first done by Hoffstein and
Rosen \cite{hr}, and later by Fisher and Friedberg \cite{ff}, whose
approach is closer to the point of view of this paper.  For $n>2$ the
$A_2$ series have been computed by Chinta \cite{a2rff}.

As stated above this theorem follows from the work of \cite{wmd1,
wmd2}.  In \cite{wmd3}, the authors study the harder problem of
constructing twisted Weyl group multiple Dirichlet series
associated to the root system $A_r$.  They construct such series for
$A_{2}$ and present a conjectural description of the series associated
to $A_r$ for arbitrary $r$ and $n$.  However, except for $\Phi = A_2$
and any $n$ or for any $A_{r}$ and $n$ sufficiently large with respect
to the rank, they are unable to prove that the series they define
satisfy the requisite functional equations.

Our method has the advantage that functional equations are essentially
built-in to our definition.  As in the case of \cite{wmd1, wmd2, wmd3,
qmds} the Weyl group multiple Dirichlet series are completely
determined by their $\p$-parts and the twisted multiplicativity
satisfied by the coefficients.  Our approach is to show that if the
$\p$-parts (which can be expressed as rational functions in the
$|\p|^{-s_i}$) satisfy certain functional equations, then the global
multiple Dirichlet series satisfies the requisite global functional
equations.  This leads us to define a certain action of the Weyl group
of the root system $\Phi$ on rational functions in $r$ indeterminates.
This approach, first introduced in \cite{chbq}, has been carried out
in the quadratic case for an arbitrary simply-laced root system, see
\cite{qmds, ratfun}.  We will extend this approach to arbitrary $\Phi$
and $n$ in forthcoming work.  However, though the basic ideas are
clear, the non-obvious group action required on rational functions can
appear unmotivated and complicated in the general setting.  Therefore,
we feel it is worthwhile in this paper to work out in detail the
simplest nontrivial case, the rank two root system $A_2.$

Here is a short plan of the paper.  Section \ref{sec:Waction}
describes the Weyl group action on rational functions that leads to a
$\p$-part \eqref{eqn:Hp} with the desired functional equations.
Although the focus of this paper is untwisted $A_{2}$, we work more
generally at first and state the full action for a general (simply
laced) root system.  We then specialize to untwisted $A_{2}$.  Section
\ref{sec:kubota} reviews the Dirichlet series of Kubota; in the
current framework, these series are Weyl group multiple Dirichlet
series attached to $A_{1}$.  The main result of this section is
Theorem \ref{thm:feE}, which shows that a certain Dirichlet series $E
(s,m)$ built from the function $H (c,d)$ \eqref{eqn:Htm} satisfies the
same functional equations as Kubota's.  Finally, in Section
\ref{sec:ddseries} we use Theorem \ref{thm:feE} to complete the proof
of Theorem \ref{thm:main}.  The basic idea is that the (one variable)
functional equations of the $E(s,m)$ induce a bivariate functional
equation in the double Dirichlet series.  

\section{A Weyl group action on rational functions}\label{sec:Waction}
Let $\Phi$ be an irreducible simply laced root system of rank $r$ with
Weyl group $W$. Choose an ordering of the roots and let $\Phi = \Phi^+
\cup \Phi^{-}$ be the decomposition into positive and negative roots.
 Let $$\Delta =\{\alpha_1, \alpha_2, \ldots, \alpha_r\}$$ be the set of
simple roots and let $\sigma_i$ be the Weyl group element
corresponding to the reflection through the hyperplane perpendicular
to $\alpha_i$.  We say that $i$ and $j$ are \emph{adjacent} if $i\neq
j$ and $(\sigma_i \sigma_j)^3=1$.  The Weyl group $W$ is generated by
the simple reflections $\sigma_1,\sigma_2,\ldots, \sigma_r$, which
satisfy the relations
\begin{equation}\label{eqn:WRelations}
(\sigma_i\sigma_j)^{r(i,j)} = 1
\mbox {\ with\ } r(i,j)=\left\{\begin{array}{ll}
3 & \mbox{\ if $i$ and $j$ are adjacent,}\\
1 & \mbox{\ if $i=j$, and}\\
2 & \mbox{\ otherwise,}
\end{array}\right.
\end{equation}
for $1\leq i, j\leq r$.
The action of the generators $\sigma_i$ on the roots is
\begin{equation}\label{eqn:WactionRoots}
\sigma_i\alpha_j=\left\{\begin{array}{ll}
\alpha_i+\alpha_j & \mbox{\ if $i$ and $j$ are adjacent,}\\
-\alpha_j & \mbox{\ if $i=j$, and  }\\
\alpha_j & \mbox{\ otherwise.}
\end{array}\right.
\end{equation}

Define
$$\sgn(w)=(-1)^{length(w)}$$
where the length function on $W$ is with respect to the generators
 $\sigma_1, \sigma_2,\ldots, \sigma_r.$
Let $\Lambda_{\Phi}$ be the lattice generated by the roots.  Any
$\alpha \in \Lambda_{\Phi}$ has a unique representation as an
integral linear combination of the simple roots:
\begin{equation}\label{eqn:repn}
\alpha=k_1\alpha_1+k_2\alpha_2+\cdots+k_r\alpha_r.
\end{equation}
We call the set $\Supp (\alpha)$ of $j$ such that $k_{j} \not = 0$ in
\eqref{eqn:repn} the \emph{support} of $\alpha$.
We denote by 
$$d(\alpha)=k_1+k_2+\cdots+k_r$$ the usual height function on $\Lambda_{\Phi}$
and put 
$$d_j(\alpha)=\sum_{i\sim j} k_i,$$
where $i\sim j$ means that the nodes labeled by $i$ and $j$ are
adjacent in the Dynkin diagram of $\Phi.$
Introduce a partial ordering on $\Lambda_{\Phi}$ by defining
$\alpha\succeq 0$ if each $k_i\geq 0$ in \eqref{eqn:repn}.  Given
$\alpha, \beta\in \Lambda_{\Phi}$, define
$\alpha\succeq\beta$ if $\alpha-\beta\succeq 0$.

Let $F=\C({\bf x})=\C(x_1,x_2,\ldots, x_r)$ be the field of rational
functions in the variables $x_1,x_2,\ldots, x_r$.  For any $\alpha \in
\Lambda_{\Phi}$, let $\x^\alpha\in F$ be the monomial
$x_1^{k_1}x_2^{k_2}\cdots x_r^{k_r}$, where the exponents $k_{i}$ are
determined as in \eqref{eqn:repn}.  Let $\ell=(l_1, \ldots , l_r) $ be
an $r$-tuple of nonnegative integers.  The tuple $\ell$ is called a
\emph{twisting parameter}; it should be thought of as corresponding to
the weight $\sum (l_{j}+1)\varpi_{j}$, where the $\varpi_{j}$ are the
fundamental weights of $\Phi$.  The case $\ell = (0,\dotsc ,0)$ is
called the \emph{untwisted case}.  For each choice of $\ell$ we will
define an action of the Weyl group $W$ on $F$.

Let $\p$ be a prime in $\OO$ of norm $p$. The action will involve
the Gauss sums $g_i(1,\p)$.\footnote{We remark that our normalization
for Gauss sums follows \cite{wmd2,wmd3} and not \cite{qmds,ratfun}.
See \cite[Remark 3.12]{qmds} for a discussion of this.}  First,
for $\x=(x_1,x_2,\ldots,
x_r)$ define $\sigma_i\x=\x^\prime$, where
\begin{equation}\label{eqn:wiaction1}
x_j^\prime=\left\{\begin{array}{ll}
px_i x_j & \mbox{\ if $i$ and $j$ are adjacent,}\\
1/(p^2 x_j) & \mbox{\ if $i=j$, and}\\
x_j & \mbox{\ otherwise.}
\end{array}\right.
\end{equation}
Next,
given $f\in F$ write 
$$ f(\x)=\sum_{\beta} a_\beta \x^\beta , $$ 
and for integers $k,i,j$ define 
\begin{equation*}
f_k(\x;i,j)=\sum_{\substack{\beta_k = i \bmod n\\
d_k(\beta) = j\bmod n}} a_\beta \x^\beta .
\end{equation*}
Finally we can define the action of $W$ on $F$ for a generator
$\sigma_k\in W$:
\begin{multline}
(f|_\ell\sigma_k)(\x)=\\
(px_k)^{l_k}\sum_{i=0}^{n-1}\sum_{j=0}^{n-1}
\left(
\PP_{ij}(x_k)f_k(\sigma_k\x;i,j-l_k)+
\QQ_{ij}(x_k)f_k(\sigma_k\x;j+1-i, j-l_k)
\right)
\end{multline}
where
\begin{eqnarray*}
  \label{eq:1}
  \PP_{ij}(x)&=&  (p x)^{1-(-2i+j+1)_n}
\frac{1-1/p}{1-p^{n-1} x^n},  \\
\QQ_{ij}(x)&=&-g^*_{2i-j-1}(1,\p)(p x)^{1-n}\frac{ 1-p^n
    x^n}{1-p^{n-1} x^n},  \\
g^*_i(1,\p) &=&\left\{
\begin{array}{ll}
  g_i(1,\p)/p & \text{if $i$ is not congruent to 0 mod $n$, }\\
-1 & \text {otherwise.}
\end{array}\right. 
\end{eqnarray*}
Here $(i)_n=i-n\lfloor i/n \rfloor$ is the minimal nonnegative integer
congruent to $i$ modulo $n$.  One can show that this action of the
generators $\sigma_{i}$ extends to
an action of $W$ on $F$; in particular the defining relations
\eqref{eqn:WRelations} are satisfied.

Now we specialize to the focus of this paper: we set $\Phi=A_2$ and
$\ell =(0,0)$.  To simplify notation we write $x, y$ for the variables
of $F$.  With these simplifications the action of $\sigma_{1}$ takes
the form
\begin{multline}
(f|\sigma_1)(x,y)=\\
\sum_{i=0}^{n-1}\sum_{j=0}^{n-1}
\left(
\PP_{ij}(x)f_1\left(\smfrac{1}{p^2x}, pxy;i,j\right)+
\QQ_{ij}(x)f_1\left(\smfrac{1}{p^2x}, pxy;j+1-i, j\right)
\right);
\end{multline}
the action of $\sigma_2$ is similar.
An invariant rational function for this action is
\begin{equation}\label{eqn:rf}
h(x,y)=\frac{N(x,y) }{(1-p^{n-1}x^n)(1-p^{n-1}y^n)
(1-p^{2n-1}x^ny^n)},
\end{equation}
where the numerator $N(x,y)$ is 
\begin{multline}\label{eqn:Np}
N (x,y) = N^{(\p)}(x,y) = 
1+g_1(1, \p) x+ g_1(1, \p) y+
pg_1(1, \p) g_2(1, \p)xy^2\\+
pg_1(1,\p)g_2(1, \p) x^2 y+ pg_1(1,  \p)^2g_2(1, \p) x^2 y^2.
\end{multline}
Let us write $h(x,y)$ as 
\begin{align}
  \label{eq:2}
  h(x,y)&=\sum_{k,l\geq 0}a(\p^k,\p^l)x^k y^l \\ \nonumber
&=\sum_{l\geq 0} y^l \left( \sum_{i=0}^{n-1}
\sum_{k = i \bmod n} a(\p^k,\p^l) x^k\right)\\ \nonumber
&=\sum_{l\geq 0}  \sum_{i=0}^{n-1} y^l \hpl(x;i),
\end{align}
say.

The following two lemmas are proved by a direct computation.
\begin{lemma}\label{lemma:hpl}
We have $N^{(\p)}(x,0)=1+g_1(1, \p) x, N^{(\p)}(0,y)=1+g_1(1, \p) y$ and
for $j = l \bmod n,$ and $0\leq i\leq n-1$,  
  \begin{equation*}
    \hpl(x;i)=(px)^l P_{ij}(x)h\left(\smfrac{1}{p^2x};i \right) 
+(px)^l Q_{ij}(x)h\left(\smfrac{1}{p^2x};l+1-i \right). 
  \end{equation*}
\end{lemma}

\begin{lemma}\label{lemma:fpl}
  Let 
$$\fpl(x;i)=\hpl(x;i)-\delta 
g_{2i-l-1}(1,\p)p^{(2i-l-2)_n}x^{(2i-l-1)_n} \hpl(x, l+1-i)$$
where $\delta=0$ if $l-2i = -1 \bmod n $ and is 1 otherwise.
Then 
$$\fpl(x;i)=(px)^{l-(l-2i)_n} \fpl\left(\smfrac{1}{p^2x};i\right).$$ 
\end{lemma}

\section{Kubota's Dirichlet series}\label{sec:kubota} The basic
building blocks of the multiple Dirichlet series are the Kubota
Dirichlet series constructed from Gauss sums \cite{kub1, kub2}.  Let
$m$ be a nonzero polynomial in $\OO$ and let $s$ be a complex
variable. These series are defined by
\begin{equation}
D(s,m)=(1-q^{n-ns})^{-1}\sum_{d\in\OOmon}\frac{g(m,d)}{|d|^s}
\end{equation}
and 
\begin{equation}
D(s,m;i)=(1-q^{n-ns})^{-1}\sum_{\substack{\deg d = i \bmod n\\ d \in \OOmon}}
\frac{g(m,d)}{|d|^s}.
\end{equation}

Kubota proved that these series have meromorphic continuation to $s\in
\C$ with possible poles only at $s=1\pm 1/n$ and satisfy a functional
equation. Actually, Kubota worked over a number field, but the
constructions over a function field are identical.
 
If the degree of $m$ is $nk+j$, where $ 0\leq j\leq n-1$, this
functional equation takes the form
\begin{equation}\label{eqn:feD}
D(s,m)=|m|^{1-s}\sum_{0\leq i\leq n-1} T_{ij}(s)D(2-s, m;i),
\end{equation}
where the $T_{ij}(s)$ are certain quotients of Dirichlet polynomials.
For fixed $s$ the $T_{ij}$ depend only on $2i-j.$ We will not need to
know anything more about the functional equation, but a more explicit
description can be found in Hoffstein \cite{jhoff} or Patterson
\cite{sjp}.

Given a set of primes $S,$ we define
\begin{equation}
D_S(s,m)=(1-q^{n-ns})^{-1}\sum_{\substack{(d,S)=1\\
d\in \OOmon }}\frac{g(m,d)}{|d|^s}.
\end{equation}
If $m_0=\prod_{\p\in S} \p$ we sometimes write $D_{m_0}(s, m)$ 
for $D_S(s,m).$

We record some properties of Gauss sums that we will use repeatedly.

\begin{prp}  Let $a,m,c,c'\in \OO.$
  \begin{list}{}{}
  \item{(i)} If $(a,c)=1$ then $g_i(am,c)=\prs{a}{c}^{-1}g_i(m,c).$
\item{(ii)} If $(c,c')=1$ then
$$g_i(m,cc')=g_i(m,c)g_i(m,c')\prs{c}{c'}^{2i}.$$
  \end{list}
\end{prp}

Using this proposition we can relate the functions $D_S$ to the
functions $D_{S'}$ for different sets $S$ and $S'.$  This is the
content of the following two lemmas.

\begin{lemma} Let $\p\in \OOmon$ be prime of norm $p.$ For an integer
$i$ with $0\leq i\leq n-1$ and $m_1,m_2,\p$ all pairwise relatively
prime, we have
$$D_{m_1}(s,m_2\p^i)=D_{\p m_1}(s, m_2\p^i)+ \frac{g(m_2\p^i,
  \p^{i+1})} 
{p^{(i+1)s}} D_{\p m_1}(s,m_2\p^{(n-i-2)_n}).$$  
More generally,
$$D(s,m)=\sum_{S_0\subset S} \Biggl (\prod_{\p\in S_0}
  \frac{g(m,\p^{i+1})}{|\p|^{(i+1)s}}\Biggr) D_S\Biggl(s, 
\prod_{\substack{\p\in S_0^c\\ \p^i||m}}\p^i\cdot
\prod_{\substack{\p\in 
  S_0\\ \p^{i}||m}}\p^{(n-i-2)_n}\Biggr). $$
\end{lemma}

\begin{proof}
We prove only the first part of the Lemma.  For $\p, m_1, m_2$ as in
the statement, 
\begin{align*}
 (1-q^{n-ns}) D_{m_1}(s,m_2\p^i)&= \sum_{\substack{(d,m_1)=1\\
d\in \OOmon }}\frac{g(m_2\p^i,d)}{|d|^s} \\
&=\sum_{k\geq 0} \sum_{\substack{(d,m_1\p)=1\\d\in \OOmon }}
\frac{g(m_2\p^i,d\p^k)}{|d|^sp^{ks}}\\
&=\sum_{k\geq 0} \sum_{\substack{(d,m_1\p)=1\\d\in \OOmon }}
\frac{g(m_2\p^i,d)g(m_2\p^i, p^k)}{|d|^sp^{ks}}\prs{d}{\p^{2k}}\\
&= \sum_{\substack{(d,m_1\p)=1\\d\in \OOmon }}
\frac{g(m_2\p^i,d)}{|d|^s}\left( \sum_{k\geq 0}
\frac{g(m_2\p^i, p^k)}{p^{ks}} \prs{d}{\p^{2k}}\right) .\\
\end{align*}
The Gauss sum in the inner sum vanishes unless $k=0$ or $i+1.$  This
proves the Lemma.
\end{proof}

Inverting the previous Lemma, we obtain

\begin{lemma}\label{lemma:dp}
If $0\leq i\leq n-2$ and $m_1, m_2, \p$ as above, 
$$D_{\p m_1}(s,m_2 \p^i)=\frac{D_{m_1}(s,m_2\p^i)}{1-|\p|^{n-1-ns}}- 
\frac{g(m_2\p^i,\p^{i+1})}{|\p|^{(i+1)s}} 
\frac{D_{m_1}(s,m_2\p^{n-i-2})}{1-|\p|^{n-1-ns}}, $$
and if $i=n-1$,
$$D_{\p m_1}(s, m_2\p^i)=
\frac{D_{m_1}(s,m_2\p^i)}{1-|\p|^{n-1-ns}}.$$ 
\end{lemma}

Now suppose that $N(x,y)=N^{(\p)}(x,y)$ is the polynomial from
(\ref{eqn:Np}).  We define a function $H$ on pairs of powers of $\p$
by setting $H(\p^k,\p^l)$ to be the coefficient of $x^{k}y^{l}$ in $N (x,y)$:
\[
N(x,y)=\sum_{}H(\p^k,\p^l)x^ky^l.
\]
We extend $H$ to all pairs of
monic polynomials by the twisted multiplicativity relation: if
$\gcd(cd,c'd')=1$, then we put
\begin{equation}
  \label{eq:16}
  H(cc',dd')=H(c,d)H(c',d')\prs{c}{c'}^2\prs{d}{d'}^2\prs{c}{d'}^{-1}
  \prs{c'}{d}^{-1}.
\end{equation}
In particular, note that
\begin{equation}
  \label{eq:19}
  H(d,1)=g(1, d).
\end{equation}

Now consider the Dirichlet series
\begin{equation}
E(s,m)=(1-q^{n-ns})^{-1}\sum_{d\in \OOmon } \frac{H(d,m)}{d^s}.
\end{equation}
That $E (s,m)$ satisfies the same functional equation as $D(s,m)$  
is the main result of this section:

\begin{theorem}\label{thm:feE}
  Let $m\in\OOmon$ be a monic polynomial of degree $nk+j$, where $0\leq j\leq
  n-1$.  Then
$$E(s,m)=|m|^{1-s}\sum_{0\leq i\leq n-1} T_{ij}(s)E(2-s, m;i).$$
\end{theorem}

\begin{proof}
Before tackling the general case, we first consider 
$m=\p^l$ for a prime $\p$.
Then
\begin{eqnarray*}
E(s,\p^l)&=& (1-q^{n-ns})^{-1}\sum_{\substack{d\in \OOmon \\
(d,\p)=1}}\sum_{k\geq 0}
\frac{H(d\p^k,\p^l)}{d^s|\p|^{ks}}\\
&=&  (1-q^{n-ns})^{-1}\sum_{\substack{d\in \OOmon \\
(d,\p)=1}}\sum_{k\geq 0}\frac{H(\p^k,
  \p^l)g(1,d)}{|\p|^{ks}d^s}\prs{d}{\p^{2k-l} }, \ \text{by
  (\ref{eq:16}) and (\ref{eq:19})}\\
&=& \sum_{k\geq 0}\frac{H(\p^k,\p^l)}{|\p|^{ks}} 
D_\p(s, \p^{(l-2k)_n})\\
&=& \sum_{j=0}^{n-1} D_\p(s, \p^{(l-2j)_n})\left(\frac{1}{|\p|^{js}}
\sum_{k\geq 0}\frac{H(\p^{j+nk},\p^l)}{|\p|^{ks}} \right)\\
&=&\sum_{j=0}^{n-1} D_\p(s, \p^{(l-2j)_n}) \hpl(|\p|^{-s};j),
\end{eqnarray*}
where $ \hpl$ was introduced in (\ref{eq:2}).
Using Lemma \ref{lemma:dp}  the previous expression becomes
\begin{equation}
  \label{eq:4}\begin{split}
&  \sum_{j=0}^{n-1}  D(s, \p^{(l-2j)_n}) \hpl(|\p|^{-s};j)\\
&- 
 \sum_{j=0}^{n-1} \delta_j \frac{g(\p^{(l-2j)_n},\p^{(l-2j)_n+1})
 }{|\p|^{((l-2j)_n+1)s}} 
 D(s, \p^{(2j-l-2)_n}) \hpl(|\p|^{-s};j),
\end{split}
\end{equation}
where $\delta_j=0$ if $l-2j\equiv n-1 (n) $ and is 1 otherwise.
Replace $j$ by $l+1-j$ in the second summation and regroup to conclude
\begin{equation}
  \label{eq:3}
  E(s,\p^l)=\sum_{j=0}^{n-1}D(s, \p^{(l-2j)_n}) \fpl(|\p|^{-s};j).
  \end{equation}
(Note the use of the identity $n-2-(l-2j)_n=(2j-l-2)_n.$)
Using the functional equations 
(\ref{eqn:feD}) of $D$ and $\fpl$ (Lemma~\ref{lemma:fpl}), we write
\begin{equation}
  \label{eq:5}\begin{split}
&   E(s,\p^l)|\p|^{-(1-s)l}\\
&=
\sum_{j=0}^{n-1}\sum_{i=0}^{n-1} T_{i, (l-2j)_n \deg \p}(s)
   D(2-s,\p^{(l-2j)_n};i)\fpl(2-s;j)\\
&= \sum_{i,j=0}^{n-1} T_{i-j\deg \p, 
(l-2j)_n \deg \p}(s)
   D(2-s,\p^{(l-2j)_n};i-j\deg \p)\fpl(2-s;j)\\
&= \sum_{i=0}^{n-1}T_{i,l\deg \p}(s)\left[ \sum_{j=0}^{n-1} 
   D(2-s,\p^{(l-2j)_n};i-j\deg \p)\fpl(2-s;j) \right]\\
&= \sum_{i=0}^{n-1}T_{i,l\deg \p}(s)
E(2-s,\p^l;i),
\end{split}
\end{equation}
where the third equality comes from our remark that the $T_{ij}$
depend only on $2i-j.$  This is the functional equation we wished to
prove, in the special case $m=\p^l.$

The argument for general $m$ is similar.  Let
$m=\p_1^{l_1}\p_2^{l_2}\cdots \p_r^{l_r}$ where the $\p_i$ are distinct
primes. Then
\begin{equation}
  \label{eq:6}
  \begin{split}
    E(s;m)&=(1-q^{n-ns})^{-1}\sum_{d\in \OOmon }\frac{ H(d,m)}{|d|^s}\\
&=(1-q^{n-ns})^{-1}\sum_{\substack{d\in \OOmon \\
(d,m)=1}} \sum_{k_1, \dots, k_r\geq 0}
\frac{ H(d\p_1^{k_1}\cdots \p_r^{k_r},\p_1^{l_1}\cdots \p_r^{l_r})}
{|d|^s|\p_1|^{k_1s}\cdots |\p_r|^{k_rs}}\\
&= (1-q^{n-ns})^{-1} \sum_{\substack{d\in \OOmon \\
(d,m)=1}} \sum_{k_1, \dots, k_r\geq 0}
\frac{ H(d,1)H(\p_1^{k_1}, \p_1^{l_1})\cdots H(\p_r^{k_r}, \p_r^{l_r})}
{|d|^s|\p_1|^{k_1s}\cdots |\p_r|^{k_rs}}  \\ 
& \qquad \times\prs{d}{m}^{-1}
\prs{d}{\p_1^{k_1}\cdots \p_r^{k_r}}^2 \prod_{a\neq b} \prs{\p_a^{k_a}}
{\p_b^{k_b}}\prs{\p_a^{l_a}} 
{\p_b^{l_b}}\prs{\p_a^{k_a}} {\p_b^{l_b}}^{-1}\\
&= \prod_{a\neq b} \prs{\p_a^{l_a}} 
{\p_b^{l_b}} \sum_{j_1=0}^{n-1}\cdots \sum_{j_r=0}^{n-1} D_m(s,
\p_1^{(l_1-2j_1)_n} \cdots \p_r^{(l_r-2j_r)_n}) \\
&\qquad \times h^{(\p_1,l_1)}(s;j_1)
\cdots h^{(\p_r,l_r)}(s;j_r)  \prod_{a\neq b} \prs{\p_a^{j_a}}
{\p_b^{j_b}}\prs{\p_a^{j_a}} {\p_b^{l_b}}^{-1}.
  \end{split}
\end{equation}
Denote for the moment by $C(j_1)=C(j_1, \ldots, j_r)$
  the product of residue symbols 
  \begin{equation}
    \label{eq:12}
    C(j_1)= \prod_{a\neq b} \prs{\p_a^{j_a}}
{\p_b^{j_b}}\prs{\p_a^{j_a}} {\p_b^{l_b}}^{-1}.
  \end{equation}

Letting $J_i=(l_i-2j_i)_n$ for $i=1, \ldots r,$ we have
\begin{equation}
  \label{eq:7}
  \begin{split}
(1-|\p_1|^{n-1-ns})&
 D_m(s,\p_1^{J_1} \cdots \p_r^{J_r})C(j_1)
= D_{\p_2\cdots \p_r}(s,\p_1^{J_1}
 \cdots \p_r^{J_r})C(j_1) \\
-&\delta_{j_1} \frac{ g(\p_1^{J_1} \cdots \p_r^{J_r}, 
\p_1^{J_1+1})}{|\p_1|^{(J_1+1)s}} 
D_{\p_2\cdots \p_r}(s,\p_1^{(2j_1-l_1-2)_n} \p_2^{J_2} 
\cdots \p_r^{J_r})C(j_1)
\end{split}
\end{equation}
by Lemma \ref{lemma:dp}.
In the second term on the right hand side, replace $j_1$ by 
$l_1+1-j_1.$  For $\delta_{j_1}\neq 0$ this gives
\begin{equation}
  \label{eq:8}
  \begin{split}
   \frac{ g(\p_1^{(2j_1-l_1-2)_n}\p_2^{J_2} 
\cdots \p_r^{J_r}, \p_1^{(2j_1-l_1-1)_n})}{|\p_1|^{((2j_1-l_1-1)_n)s}} 
D_{\p_2\cdots \p_r}(s,\p_1^{J_1} \p_2^{J_2} 
\cdots \p_r^{J_r})C(l_1-j_1+1).
  \end{split}
\end{equation}
The Gauss sum can be written as 
\begin{equation}
  \label{eq:9}
  \prs{\p_2^{J_2}\cdots \p_r^{J_r}}{\p_1^{2j_1-l_1-1}}^{-1}
g(\p_1^{(2j_1-l_1-2)_n}, \p_1^{(2j_1-l_1-1)_n}),
\end{equation}
and $C(l_1-j_1+1)$ is 
\begin{equation}
  \label{eq:10}
  \prs{\p_2^{J_2}\cdots \p_r^{J_r}}{\p_1^{l_1-j_1+1}}^{-1}
 \Biggl(\prod_{\substack{a\neq b\\ a,b\neq 1}} \prs{\p_a^{j_a}}
{\p_b^{j_b}}\prs{\p_a^{j_a}} {\p_b^{l_b}}^{-1}\Biggr).
  \end{equation}
Taking the product of (\ref{eq:9}) and (\ref{eq:10}) yields
\begin{equation}
  \label{eq:11}\begin{split}
   \prs{\p_2^{J_2}\cdots \p_r^{J_r}}{\p_1^{j_1}}^{-1}  &
g(\p_1^{(2j_1-l_1-2)_n}, \p_1^{(2j_1-l_1-1)_n}) \Biggl(\prod_{\substack{a\neq b\\ a,b\neq 1}} \prs{\p_a^{j_a}}
{\p_b^{j_b}}\prs{\p_a^{j_a}} {\p_b^{l_b}}^{-1}\Biggr)\\
  &= g(\p_1^{(2j_1-l_1-2)_n}, \p_1^{(2j_1-l_1-1)_n})C(j_1).
\end{split}\end{equation}
Therefore, continuing from the last line of (\ref{eq:6}),
\begin{equation}
  \label{eq:13}
  \begin{split}
    E(s,m)&=\prod_{a\neq b} \prs{\p_a^{l_a}} 
{\p_b^{l_b}} \sum_{j_1=0}^{n-1}\cdots \sum_{j_r=0}^{n-1} D_{m'}(s,
\p_1^{(l_1-2j_1)_n} \cdots \p_r^{(l_r-2j_r)_n}) \\
&\times \prod_{a\neq b} \prs{\p_a^{j_a}}
{\p_b^{j_b}}\prs{\p_a^{j_a}} {\p_b^{l_b}}^{-1}
 f^{(\p_1,l_1)}(s;j_1) h^{(\p_2,l_2)}(s;j_2)
\cdots h^{(\p_r,l_r)}(s;j_r),
  \end{split}
\end{equation}
where $m'=\p_2^{l_2}\cdots \p_r^{l_r}.$  Repeating this procedure to
remove the primes from $m$ one at a time, we find that up to a
constant of modulus one, $E(s,m)$ is equal to 
\begin{equation}
  \label{eq:14}
    \begin{split}
   \sum_{j_1=0}^{n-1}\cdots \sum_{j_r=0}^{n-1} D_{}(s,
\p_1^{J_1} \cdots \p_r^{J_r})\left( 
\prod_{a=1}^{r} f^{(\p_a,l_a)}(s;j_a)  \right)
\prod_{a\neq b} \prs{\p_a^{j_a}}
{\p_b^{j_b}}\prs{\p_a^{j_a}} {\p_b^{l_b}}^{-1}.
   \end{split}
\end{equation}
We may now apply the functional equations of $D$ and the
$f^{(\p_a,l_a)}$ as in (\ref{eq:5}) to 
conclude that $E(s, m)$ satisfies the functional
equation 
\begin{equation}
  \label{eq:15}
  E(s,m)= |m|^{1-s}\sum_{i=0}^{n-1}T_{i,\deg m}(s)
E(2-s,m;i).
\end{equation}
This completes the proof of the theorem.
\end{proof}

For later use, we record the following bound:
\begin{prp}
  \label{prop:Ebound}
For all $\epsilon>0,$   $m\in\OO$ and $0\leq i <n,$
\begin{equation*}
 (s-1-\smfrac 1n)(s-1+\smfrac 1n)
 E(s, m;i) \ll_\epsilon \left\{
    \begin{array}{ll}
      1 & \text{for $\re(s)>\smfrac 32+\epsilon$}\\
      |m|^{\smfrac 12+\epsilon} & 
\text{for $\smfrac 12-\epsilon<\re(s)<\smfrac 32+\epsilon$}\\
|m|^{1-s+\epsilon} & \text{for $\re(s)<\smfrac 12-\epsilon$}\\
    \end{array}
\right.
\end{equation*}
\end{prp}

\begin{proof}
  Use the meromorphy and functional equation of $E(s,m)$ together with
  the convexity principle, cf.  \cite[Eq.~(2.3)]{fhl} and
  \cite[Propostion 8.4]{lp}.
\end{proof}

\section{The double dirichlet series}\label{sec:ddseries}

Recall the definition of the double Dirichlet series from
\eqref{def:Z}--\eqref{eqn:Hp}.  In this section we show that $Z(s_1,s_2)$
has a meromorphic continuation to $s_1,s_2\in \C$ and satisfies a
group of functional equation isomorphic to $W.$  In \cite{wmd1}, the
authors show in detail how the analytic continuation of a Weyl group
multiple Dirichlet series follows from the functional
equations. Therefore we concentrate on establishing the functional
equations of $Z(s_1, s_2).$ 

Actually we need to consider slightly
different series.  For integers $0\leq i,j\leq n-1$ we define
\begin{equation}
  \label{eq:17}\begin{split}
  &Z(s_1,s_2;i,j)=\\ &(1-q^{n-ns_1})^{-1}\! (1-q^{n-ns_2})^{-1}\!
(1-q^{2n-ns_1-ns_2})^{-1}\!\!\! \!\!\!
\sum_{\substack{m\in \OOmon \\
 \deg m = i \bmod n}} \sum_{\substack{d\in \OOmon \\
 \deg d = j \bmod n}}\!\!
\frac{H(d,m)}{|m|^{s_1}|d|^{s_2}}.
\end{split}
\end{equation}
We further introduce the notation
\begin{equation*}
  Z(s_1, s_2;i, * )=\sum_{j} Z(s_1, s_2;i, j)
\end{equation*}
and 
\begin{equation*}
  Z(s_1, s_2;*, j )=\sum_{i} Z(s_1, s_2;i, j).
\end{equation*}
These series are absolutely convergent for $\re(s_1), \re(s_2)>3/2.$
In fact, we can do a little better.   Summing over $d$ first yields
 \begin{align}\nonumber
  Z(s_1, s_2;i,*)&=(1-q^{n-ns_1})^{-1}\! (1-q^{n-ns_2})^{-1}\!
(1-q^{2n-ns_1-ns_2})^{-1}\\ \nonumber
&\qquad\times \sum_{\substack{m\in \OOmon \\
 \deg m = i \bmod n}}
\left ( \frac{1}{|m|^{s_1}}\sum_{\substack{d\in \OOmon \\}}
 \frac{H(d,m)}{|d|^{s_2}} \right) \\   \label{eq:18}
&=(1-q^{n-ns_1})^{-1}\! (1-q^{2n-ns_1-ns_2})^{-1}
\sum_{\substack{m\in \OOmon \\
 \deg m = i \bmod n}}
\frac{E(s_2, m)}{|m|^{s_1}}
\end{align}
By the convexity bound of Proposition \ref{prop:Ebound}, this representation
of $Z(s_1,s_2;i, *)$ is seen to meromorphic for $\re(s_1)>0, \re
(s_2)>2.$   Alternatively, summing over $m$ first we deduce that
$Z(s_1,s_2;i, *)$ is meromorphic for $\re(s_2)>0, \re (s_1)>2.$ 
Let $\mathcal R$ be the tube domain that is the union of these three
regions of initial meromorphy:  
\begin{multline*}
  {\mathcal R}=\{\re(s_1), \re(s_2)>3/2\} \cup \{\re(s_1)>0, \re
(s_2)>2\} \\
\cup \{\re(s_2)>0, \re (s_1)>2\}.
\end{multline*}
Let the Weyl group $W$ act on $\C^2$ by 
\begin{equation}\label{eq:WonC2}
  \sigma_1:(s_1, s_2)\mapsto (2-s_1, s_1+s_2-1), \ \ \sigma_2:(s_1,
  s_2)\mapsto (s_1+s_2-1, 2-s_2).
\end{equation}
Let $\mathcal F$ be the real points of a closed fundamental
domain for the action of $W$ on $\C^{2}$: 
\[
\mathcal{F} = \{\re (s_{1}), \re (s_{2})\geq 1 \}.
\]
One can easily see that
$\mathcal R \smallsetminus \mathcal F \cap \mathcal R$ is compact.
Therefore, by the principle of analytic continuation and Bochner's
tube theorem \cite{bo}, to prove that $Z(s_1, s_2)$ has a meromorphic
continuation to $\C^2$ it suffices to show that the functions
$Z(s_1,s_2;i,j)$ satisfy functional equations as $(s_1, s_2)$ goes to
$ (2-s_1, s_1+s_2-1)$ and $ ( s_1+s_2-1, 2-s_2).$ For details, we
refer to \cite [Section 3]{wmd1}.

To prove the $\sigma_2$ functional equation, we begin with
(\ref{eq:18}) and write
\begin{align*}
  Z(s_1, s_2;i,*)&=(1-q^{n-ns_1})^{-1}\! (1-q^{2n-ns_1-ns_2})^{-1}
\sum_{\substack{m\in \OOmon \\
 \deg m = i \bmod n}}
\frac{E(s_2, m)}{|m|^{s_1}}\\
&=(1-q^{n-ns_1})^{-1}(1-q^{2n-ns_1-ns_2})^{-1} \\
&\qquad\times\!\!\!\sum_{\substack{m\in \OOmon \\
 \deg m = i \bmod n}}
\!\!\!\frac{|m|^{1-s_2}}{|m|^{s_1}}\sum_{j=0}^{n-1}T_{ji}(s_2)
E(2-s_2,m;j), \ \ \text{by Thm. \ref{thm:feE}}\\
&=\sum_{j=0}^{n-1}T_{ji}(s_2) Z(s_1+s_2-1, 2-s_2;i,j)
\end{align*}
The $\sigma_1$ functional equation is proved similarly.

We conclude that

\begin{theorem}\label{thm:final}
The double Dirichlet series has a meromorphic continuation to
$s_1,s_2\in \C$ and is holomorphic away from the hyperplanes
\begin{equation*}
  s_1=1\pm \smfrac 1n, s_2=1\pm \smfrac 1n \text{\ and\ }
s_1+s_2=2\pm \smfrac 1n. 
\end{equation*}
Furthermore, $Z(s_1, s_2)$ satisfies the functional equations
\begin{align*}
  Z(s_1, s_2)& =\sum_{i,j}T_{ji}(s_2) Z(s_1+s_2-1, 2-s_2;i,j) \\
& =\sum_{i,j}T_{ij}(s_1) Z(2-s_1,s_1+s_2-1;i,j).
\end{align*}
\end{theorem}

\bibliographystyle{amsplain}

\end {document}